\documentclass[12pt]{article}
\usepackage{amssymb}

\input{tcilatex}
\begin{document}

\bigskip \ 

\bigskip \ 

\begin{center}
\textbf{SOME EXCEPTIONAL CASES IN MATHEMATICS: }

\textbf{\smallskip \ }

\textbf{EULER CHARACTERISTIC, DIVISION ALGEBRAS, }

\smallskip \ 

\textbf{CROSS VECTOR PRODUCT AND FANO MATROID}

\smallskip \ 

J. A. Nieto\footnote{%
nieto@uas.uasnet.mx, janieto1@asu.edu}

\smallskip

\textit{Mathematical, Computational \& Modeling Science Center, Arizona
State University, PO Box 871904, Tempe, AZ 85287, USA}

\smallskip \ 

\textit{Facultad de Ciencias F\'{\i}sico-Matem\'{a}ticas de la Universidad
Aut\'{o}noma de Sinaloa, 80010, Culiac\'{a}n Sinaloa, M\'{e}xico}

\smallskip \ 

\textit{Departamento de Investigaci\'{o}n en F\'{\i}sica de la Universidad
de Sonora, 83000, Hermosillo Sonora , M\'{e}xico}

\bigskip \ 

Abstract

\smallskip \ 
\end{center}

We review remarkable results in several mathematical scenarios, including
graph theory, division algebras, cross product formalism and matroid theory.
Specifically, we mention the following subjects: (1) the Euler relation in
graph theory, and its higher-dimensional generalization, (2) the dimensional
theorem for division algebras and in particular the Hurwitz theorem for
normed division algebras, (3) the vector cross product dimensional
possibilities, (4) some theorems for graphs and matroids. Our main goal is
to motivate a possible research work in these four topics, putting special
interest in their possible links.

\bigskip \ 


November, 2011

\newpage

\noindent \textbf{1. Introduction}

\smallskip \ 

There is not doubt that in mathematics exceptional cases are always very
attractive subjects. As Stillwell[1] has remarked \textquotedblleft in the
mind of every mathematician, there is a tension between the general rule and
exceptional cases. Our conscience tell us we should strive for general
theorems, yet we are fascinated and seduced by beautiful
exceptions\textquotedblright . He adds that the solution of this dilemma
could be a general theory of exceptions which contains a complete
description of their structure and relations. However, such a general theory
seems far from reaching. At present it seems more reasonable to look for a
subset of all possible exceptions which we feel may have some kind of
relation. In this sense, in this work, we shall expose four remarkable
scenarios, each one in apparently unrelated areas. We refer to the following
scenarios: (1) the Euler relation (see Ref. [2] and references therein) in
graph theory [3], and its higher-dimensional generalization [4], (2) The
Milnor-Bott-Kervaire [5], [6] theorem for division algebras and the Hurwitz
theorem [7] for normed division algebras (see also Refs. [8] and [9]), (3)
the vector cross product dimensional possibilities (see Ref. [10] and
references therein) and (4) some interesting theorems for graphs and
matroids [11], including two Whitney theorems [12]-[13] (see also Ref.
[14]). The main idea is to motivate a research work on a possible relation
between these four topics. In fact, since these subjects are exceptional one
should expect to find special motivation in a quest for links between them.
At this respect, our conjecture is that the abstract duality concept
[15]-[16] may be the key mathematical notion for finding such links.

\bigskip \ 

\noindent \textbf{2. Euler relation}

\smallskip \ 

Descartes in 1640 and Euler in 1752 (see Refs. [2] and [4]) observed a
fundamental relation between the number of vertices $V$, edges $E$ and faces 
$F$ of a polyhedra in three dimensions. Euler expressed this important
geometrical fact in the famous formula $V-E+F=2$. Since the Euler's
discovery, generations of mathematicians have been fascinated for this
result, including Poincar\'{e} himself.

At present the above result is expressed as a theorem:

\bigskip \ 

\bigskip \ 

\smallskip \ 

\noindent \textbf{I. Theorem} (\textbf{Euler, 1752}):

\smallskip \ 

Let $G$ be a connected planar graph, and let $V$, $E$ and $F$ denote
respectively the number of vertices, edges and faces of $G$. Then

\begin{equation}
\chi \equiv V-E+F=2.  \tag{1}
\end{equation}%
The symbol $\chi $ is known as the Euler characteristic.

\smallskip \ 

\textit{Proof}: There are a number of proofs available in the literature for
this theorem. One of the simplest and, perhaps one of the more elegant, is
the one that is based in the dual graph $G^{\ast }$ of a planar graph $G$
(see Ref. [12]), which is defined as follows:

\begin{equation}
F^{\ast }\equiv V,  \tag{2}
\end{equation}%
\begin{equation}
E^{\ast }\equiv E,  \tag{3}
\end{equation}%
and%
\begin{equation}
V^{\ast }\equiv F,  \tag{4}
\end{equation}%
where $V^{\ast }$,$E^{\ast }$ and $F^{\ast }$ are the vertices, edges and
faces of $G^{\ast }$ respectively. Geometrically, for the construction of $%
G^{\ast }$ one chooses a vertex $V^{\ast }$ in each face $F$ of $G$. Given
two faces $F_{1}$and $F_{2}$ of $G$, we join the corresponding vertices $%
V_{1}^{\ast }$ and $V_{2}^{\ast }$ by one and only one edge $E^{\ast }$
which crosses the edge $E$, determined by the boundary of the two faces $%
F_{1}$and $F_{2}$.

We shall define the two fundamental quantities

\begin{equation}
R=V-1  \tag{5}
\end{equation}%
and

\begin{equation}
N=E-R,  \tag{6}
\end{equation}%
which are called the rank and nullity (or corank, or first Betti number) of $%
G$, respectively. Similarly, one may consider the rank $R^{\ast }$ and
nullity $N^{\ast }$ of the dual graph $G^{\ast }$ as follows:

\begin{equation}
R^{\ast }=V^{\ast }-1  \tag{7}
\end{equation}%
and

\begin{equation}
N^{\ast }=E^{\ast }-R^{\ast }.  \tag{8}
\end{equation}%
(When $G$ is not connected we write

\begin{equation}
R=V-k,  \tag{9}
\end{equation}%
where $k$ denotes the number of connected pieces.)

By virtue of the Euler characteristic (1), which is true for connected
planar graphs we find%
\begin{equation}
R^{\ast }=N  \tag{10}
\end{equation}%
and

\begin{equation}
N^{\ast }=R.  \tag{11}
\end{equation}%
Conversely given (10) and (11) one can derive (1). In fact, one first
observes that

\begin{equation}
\chi =V-E+F=V-E^{\ast }+V^{\ast }.  \tag{12}
\end{equation}%
So, we have%
\begin{equation}
\begin{array}{c}
\chi =(V-1)-E^{\ast }+(V^{\ast }-1)+2 \\ 
\\ 
=R-N^{\ast }+2=2,%
\end{array}
\tag{13}
\end{equation}%
where in the last step we used (11).

It is interesting to observe that one can also write (1) in form%
\begin{equation}
\chi =V-E+F=F^{\ast }-E+F.  \tag{14}
\end{equation}%
So, one sees from (12) and (14) that the the Euler characteristic $\chi $
relates the vertices $V$ and the dual vertices $V^{\ast }$, or the faces $F$
and the dual faces $F^{\ast }$, with the selfdual edges$\ E=E^{\ast }$.

In the case that the planar graph is not connected the Euler characteristic
(1) is generalized in the form

\begin{equation}
\chi =V-E+F=1+k.  \tag{15}
\end{equation}

Using (1) one can prove the fact that there are only five regular
polyhedrons (the so called Platonic solids): In fact, for regular
polyhedrons we have the relations

\begin{equation}
pF=2E  \tag{16}
\end{equation}%
and

\begin{equation}
qV=2E,  \tag{17}
\end{equation}%
where $p$ is the number of edges in a regular polygon and $q$ is the number
of faces in a given vertex. Substituting (16) and (17) into (1) we find

\begin{equation}
\frac{2E}{q}-E+\frac{2E}{p}=2  \tag{18}
\end{equation}%
or

\begin{equation}
E(\frac{2p-qp+2q}{pq})=2.  \tag{19}
\end{equation}%
Since $E>0$ and $pq>0$ we have $2p-qp+2q>0$ or $-2p+qp-2q<0$ and therefore
we get the formula%
\begin{equation}
(q-2)(p-2)<4.  \tag{20}
\end{equation}%
The only possible solution of (20) are%
\begin{equation}
\begin{array}{c}
\text{(I) }(p-2)=1\text{ and }(q-2)=1,\text{ then }p=3\text{ and }q=3\text{
(tetrahedron).} \\ 
\text{(II) }(p-2)=2\text{ and }(q-2)=1,\text{ then }p=4\text{ and }q=3\text{
(cube). }\  \text{ \thinspace \qquad } \\ 
\text{(III) }(p-2)=1\text{ and }(q-2)=2,\text{ then }p=4\text{ and }q=3\text{
(octahedron). \thinspace \ } \\ 
\, \, \text{(IV) }(p-2)=3\text{ and }(q-2)=1,\text{ then }p=5\text{ and }q=3%
\text{ (dodecahedron).} \\ 
\text{(V) }(p-2)=1\text{ and }(q-2)=3,\text{ then }p=3\text{ and }q=5\text{
(icosahedron).}%
\end{array}
\tag{21}
\end{equation}

\noindent Note that (16) and (17) can also be written as

\begin{equation}
pV^{\ast }=2E^{\ast }  \tag{22}
\end{equation}%
and

\begin{equation}
qF^{\ast }=2E^{\ast }.  \tag{23}
\end{equation}%
Therefore, in the context of the dual graph $G^{\ast }$, $p$ is the number
of faces in a given dual vertex and $q$ can be interpreted as the number of
edges of a polygon. As expected, from this point of view (II) and (III) are
just dual polyhedra. Similarly, (IV) and (V) are also duals, while (I) is a
self-dual case.

A generalization to higher dimensions $n$ of the Euler characteristic $\chi $
is given by (see Ref. [17])

\begin{equation}
\chi (M)=\dsum \limits_{i=1}^{n}(-1)^{i}\alpha _{i}.  \tag{24}
\end{equation}%
Here, $\alpha _{i}$ denotes the $i$-simplex associated with a $n$%
-dimensional simplicial complex manifold $M$. One can show that (24) can
also be written in terms of the Betti numbers $b_{i}$,

\begin{equation}
\chi (M)=\dsum \limits_{i=1}^{n}(-1)^{i}b_{i}.  \tag{25}
\end{equation}%
In this case $\chi (M)$ is called Euler-Poincar\'{e} characteristic of $M$.

Moreover, one has the generalized Gauss-Bonnet result

\begin{equation}
\chi (M)=\int_{M}e.  \tag{26}
\end{equation}%
Here, we have

\begin{equation}
e=\varepsilon ^{i_{1}...i_{n}}\varepsilon
_{a_{1}...a_{n}}R_{i_{1}i_{2}}^{a_{1}a_{2}}...R_{i_{n-1}i_{n}}^{a_{n-1}a_{n}},
\tag{27}
\end{equation}%
with $R_{i_{1}i_{2}}^{a_{1}a_{2}}$ the Riemann curvature tensor and the $%
\varepsilon $-symbol $\varepsilon ^{i_{1}...i_{n}}$ is the completely
antisymmetric density tensor.

If $M=S^{2}$, where $S^{2}$ is the $2$-sphere, then we get

\begin{equation}
\chi (S^{2})=\frac{1}{2\pi }\int_{S^{2}}\varepsilon ^{ij}\varepsilon
_{ab}R_{ij}^{ab}=\frac{1}{4\pi }\int_{S^{2}}\sqrt{g}R=2.  \tag{28}
\end{equation}%
In general, for the $n$-sphere we have $\chi (S^{n})=2$ if $n$ is even or $%
\chi (S^{n})=0$ if $n$ is odd.

There is another surprising connection between an arbitrary vector field on
a surface and the Euler characteristic $\chi $. This comes with the name of
Poincar\'{e} index theorem:

\smallskip \ 

\noindent \textbf{Poincar\'{e} Index Theorem: }Let $\mathcal{V}$ be a
tangent vector field on a smooth surface $S$ with only isolated critical
points $i=1,...,k$. Then

\begin{equation}
\dsum \limits_{i=1}^{k}I_{i}=\chi (S).  \tag{29}
\end{equation}%
For $S^{2}$ we have

\begin{equation}
\dsum \limits_{i=1}^{k}I_{i}=\sin k-saddle+source=\chi (S).  \tag{30}
\end{equation}%
There is a canonical way to understand this relation. Place a sink in the
middle of each triangle, a source at each vertex and saddle point at the
midpoint of each edge.

Moreover, for a closed orientable surface $S_{g}$ of genus $g$, that is,
with $g$ holes or handles we have

\begin{equation}
\chi (S_{g})=2-2g.  \tag{31}
\end{equation}%
For the sphere one has no holes, so $g=0$ and therefore $\chi (S_{0})=2$ as
expected.

\bigskip \ 

\noindent \textbf{3. Division Algebras}

\smallskip \ 

One of the most remarkable theorems in topology is the following [5] and [6]:

\smallskip \ 

\noindent \textbf{II.} \textbf{Theorem (Milnor-Bott-Kervaire): }The only
dimensions $n$ for which we have multiplication $R^{n}\times
R^{n}\rightarrow R^{n}$ denoted with $xy=0$ implying either $x=0$ or $y=0$
are $n=1,2,4$ or $8$. These multiplication can be realized respectively by
the real numbers $R$, the complex numbers $C$, the quaternions $H$, and
octonions $O$.

\smallskip \ 

In fact, this theorem turns out to be a generalization of the Hurwitz
theorem for normed division algebras [7]:

\smallskip \ 

\noindent \textbf{Theorem (Hurwitz 1986): }Every normed division algebra
with an identity is isomorphic to one of the following algebras: $R$, $C$, $%
H $ and $O$.

\smallskip

\noindent In turn, this is closely related to the fact that the only
parallelizable spheres are $S^{1},S^{3}$ and $S^{7}$ (see Refs. [17] and
[18] and references therein). We also find a connection between the Hurwitz
theorem and the generalized Frobenius theorem, namely (see Ref. [19] and
references therein)

\smallskip \ 

\noindent \textbf{Theorem (Frobenius): }Every alternative division algebra
is isomorphic to one of the following algebras: $R$, $C$, $H$ and $O$.

\smallskip \ 

The proof of the Milnor-Bott-Kervaire theorem uses the methods of
characteristic classes [20], but the key element in the proof is the Bott
periodicity theorem [21]. An interesting observation is that Bott
periodicity is related to the parallelizable property of spheres [22]. This
has motivated to make a formal study of vectors fields on spheres [23] and
to develop the so-called K-theory [24].

\bigskip \ 

\smallskip \ 

\noindent \textbf{4. Vector Cross Product}

\smallskip \ 

Here, we simply mention the following theorem (see Ref. [10] and references
therein)

\smallskip \ 

\noindent \textbf{III. Theorem (Generalized vector cross product) }%
Alternating vector cross product $R^{n}\times R^{n}...\times
R^{n}=R^{rn}\rightarrow R^{n}$ is possible only in the following cases:

\smallskip \ 
\begin{equation}
\begin{array}{c}
\text{(a) }r=1,n\text{ even}.\qquad \qquad \, \\ 
\\ 
\text{(b) }r=n-1.\qquad \qquad \quad \quad \\ 
\\ 
\text{(c) }r=2,n=3\text{ or }n=7.\text{ }\quad \\ 
\\ 
\text{(d) }r=3,n=8.\qquad \qquad \quad%
\end{array}
\tag{32}
\end{equation}%
It is not difficult to see that the cases $r=3$, $n=8$ and $r=2$, or $n=7$
are closely related to octonions $O$, while the case $r=2$, or $n=3$ is
connected with quaternions $H$. Moreover, the case $r=1$, $n=2$ is related
with complex numbers $C$. So the case (a) corresponds to copies of complex
structure $C$. The only new aspect is the case (b) for $n\neq 3$, which can
be described in terms of the $\varepsilon $-symbol. So, the possible vector
cross products are closely related again with the existence of $R,C,H$ and $%
O $.

\smallskip \ 

\noindent \textbf{5. Matroid theory}

\smallskip \ 

Presumable, matroid theory emerges in the year 1935 with the work of Whitney
[13] on the abstract properties of linear dependence. If one compares such a
work with the Whitney's paper [12] on graph theory published in 1932, one
observes some kind of influence of the following theorem:

\smallskip \ 

\noindent \textbf{Theorem (Kuratowski, 1930); }A finite graph is planar if
and only if it does not contain a subgraph that is homeomorphic to $K_{5}$
or $K_{3,3}$.

\smallskip \ 

\noindent Here $K_{5}$ denotes the complete graph of five vertices and $%
K_{3,3}$ is the binary matroid of three vertices. A related theorem
establishes that

\smallskip \ 

\noindent \textbf{\ Theorem (Planar-Dual); }A finite graph is planar if and
only if it has a dual.

\smallskip \ 

\noindent So, a graph has a dual if and only if it does not contain a
subgraph that is homeomorphic to $K_{5}$ or $K_{3,3}$. Perhaps, Whitney
discovered the concept of a matroid by insisting in making sense of a
possible duality for $K_{5}$ and $K_{3,3}$. In fact, a matroid $M$ is a pair 
$(E,\mathcal{B})$, where $E$ is a finite set and $\mathcal{B}$ is a
collection of subsets of $E$, called "bases", with the following properties
(see Ref. [11] for details):

\smallskip \ 

(A) $\mathcal{B}$ is nonempty.

(B) No member of $\mathcal{B}$ is a proper subset of another.

(C) If $B$ and $B^{\prime }$ are distinct members of $\mathcal{B}$ and $b$
is an element of $B$ not belonging to $B^{\prime }$, then there exists an
element $b^{\prime }$ belonging to $B\backslash B^{\prime }$ such that $%
B-b\cup b^{\prime }$ is a basis. (This property is called the basis exchange
property.)

\smallskip \ 

\noindent A dual $M^{\ast }$ of $M$ is defined in terms of the pair $(E,%
\mathcal{B}^{\ast })$, where $\mathcal{B}^{\ast }=\{B^{\ast }\mid B^{\ast
}=E\backslash B\}$. From this definition one can show that every matroid $M$
has a unique dual $M^{\ast }$. So, the matroids $M(K_{5})$ and $M(K_{3,3})$
associated with the graphs $K_{5}$ and $K_{3,3}$ respectively must have a
dual $M^{\ast }(K_{5})$ and $M^{\ast }(K_{3,3}$). The surprising thing is
that now there is not graphs associated with $M^{\ast }(K_{5})$ and $M^{\ast
}(K_{3,3}$). In fact, the matroids $M^{\ast }(K_{5})$ and $M^{\ast }(K_{3,3}$%
) are not graphic. Formally one has the following interesting results:

\smallskip \ 

\noindent \textbf{\ Theorem; }Let $G$ be a graph. Then $G$ is planar if an
only if $M(G)$ has not a submatroid isomorphic to $M(K_{5})$ and $M(K_{3,3}$%
).

\smallskip \ 

\noindent \textbf{\ Theorem; }A matroid is graphic if and only if it has no
a submatroid isomorphic to any of the matroids $U_{2,4},$ $F_{7}$, $%
F_{7}^{\ast }$, $M(K_{5})$ and $M(K_{3,3})$.

\smallskip \ 

\noindent \textbf{\ }Here, $U_{2,4}$ is the uniform matroid and $F_{7}$ is
the Fano matroid (see Ref. [11] for details) defined as follows

\smallskip \ 

\noindent \textbf{\ Definition (Fano Matroid); }The Fano matroid is a
matroid defined on the set $E=\{1,2,3,4,5,6,7\}$, whose bases are all the
triples of $E$ except those determined by $\{1,2,4\}$, $\{2,3,5\}$, $%
\{3,4,6\}$, $\{4,5,7\}$, $\{5,6,1\}$, $\{6,7,2\}$ and $\{7,1,3\}$.

\smallskip \ 

\noindent \textbf{\ IV. Theorem (Fano matroid); }The Fano matroid is binary,
but not graphic, cographic, transversal, or regular.

\smallskip \ 

\noindent (See Refs. [3] and [11] for details.

\noindent So, the matroid concept can no only be understood as a
generalization of the graph concept but also can be viewed as a mathematical
structure in which the duality symmetry still plays a central role.

Two interesting theorems in graph theory that, are also true in matroid
theory, are:

\smallskip \ 

\noindent \textbf{\ Theorem (Whitney, 1932); }Let $G_{1}$, ...,$G_{m}$ and $%
G_{1}^{\ast }$, ..., $G_{m}^{\ast }$ be the blocks of $G$ and $G^{\ast }$
respectively, and $G_{i}^{\ast }$ be dual of $G_{i}$, $i=1,...,m$. Then, $G$
is dual of $G^{\ast }$.

\smallskip \ 

\noindent \textbf{Theorem (Whitney, 1932); }Let $G_{1}$, ...,$G_{m}$ and $%
G_{1}^{\ast }$, ..., $G_{m}^{\ast }$ be the blocks of the dual graphs $G$
and $G^{\ast }$ and let the correspondence between these graphs be such that
edges in $G_{i}^{\ast }$ corresponds to edges in $G_{i}$, $i=1,...,m$. Then $%
G_{i}^{\ast }$ is dual of $G_{i}$.

\smallskip \ 

\noindent At present there are a enormous amount of information in the
literature about matroid theory (see Ref. [11] and references therein). But
perhaps one of the most interesting developments is oriented matroid theory
[25], which in the context of graph theory corresponds to oriented (or
directed) graphs. (For an application of oriented matroids to high energy
physics see Refs [26]-[27] and references therein.)

\smallskip \ 

\noindent \textbf{6. Final comments}

\smallskip \ 

In view of section (1), one may wonders whether the generalization of the
Euler characteristic given (1) can also be proved using duality. In
particular, one may be interest to see if matroid theory can help to show
Euler formula for polytopes. This will be relevant because as we mentioned
before every matroid has a dual and therefore one should expect that also
duality is linked to polytopes. This has been in fact what Lawrence [28]
proved using an algebraic methods.

Similarly, one would like to know whether the formula (31) for Riemann
surfaces of genus $g$ can be proved using the duality concept. In fact, we
are aware that formula (31) has been proved by various methods, but we have
not seen a proof using duality as is the case of formula (1). In what
follows, we shall explain a possible prove of (31) using duality.

Let us first write (31) in the form%
\begin{equation}
\chi +2g=2.  \tag{33}
\end{equation}%
If one assumes that $\chi $ can be written as in (12), that is, as

\begin{equation}
\chi =V-E^{\ast }+V^{\ast },  \tag{34}
\end{equation}%
we get the formula

\begin{equation}
V-E^{\ast }+V^{\ast }+2g=2  \tag{35}
\end{equation}%
or

\begin{equation}
\mathcal{V}-E^{\ast }+\mathcal{V}^{\ast }=2,  \tag{36}
\end{equation}%
where

\begin{equation}
\begin{array}{c}
\mathcal{V}=V+g, \\ 
\\ 
\mathcal{V}^{\ast }=V^{\ast }+g.%
\end{array}
\tag{37}
\end{equation}%
So, assuming again the dualities relations

\begin{equation}
\begin{array}{c}
\mathcal{R}=\mathcal{V}-1=\mathcal{N}^{\ast }, \\ 
\\ 
\mathcal{R}^{\ast }=\mathcal{V}^{\ast }-1=\mathcal{N},%
\end{array}
\tag{38}
\end{equation}%
with

\begin{equation}
\begin{array}{c}
\mathcal{N}=E-\mathcal{R}, \\ 
\\ 
\mathcal{N}^{\ast }=E^{\ast }-\mathcal{R}^{\ast },%
\end{array}
\tag{39}
\end{equation}%
one discovers, following similar procedure as in section 2, that (33) holds.
It remains to clarify expressions (36). The first observation is that we did
not need to introduce dual faces as in (2). This is because one can draw a
graph in a surface of genus $g$ without faces. And this means that we can
not define duality vertices in the sense of expression (2). However, one may
think in associating a virtual vertices to a circuit of a graph with no
face. It turns out that the minimum possible virtual vertices for $G$ is $g$%
, and dually the minimum number of virtual vertices for $G^{\ast }$ is also $%
g$. This is the sense of (36) and (37).

Finally, we would like to point out the possibility that and abstract
duality may be the key concept for a link between the exceptional cases
discussed in sections 2-5. A general definition of abstract duality for
matroids has been given by Bland and Dietrich [15]-[16]. Let $\mathcal{M}$
denote the family of all matroids $M$ on a finite set of elements $E$
(ground set). The matroid duality relation $D:\mathcal{M\longrightarrow M}$
is an involution

\smallskip \ 

\noindent (I) $D(D(M))=M$ (for all $M\in \mathcal{M}$),

\smallskip \ 

\noindent that preserves the ground set

\smallskip \ 

\noindent (II) $E(D(M))=E(M)$ (for all $M\in \mathcal{M}$),

\smallskip \ 

\noindent and

\smallskip \ 

\noindent (III) (a) $D(M\backslash e)=D(M)/e$ (for all $M\in \mathcal{M}$, $%
e\in E$).

\thinspace \thinspace \thinspace \ (b) $D(M/e)=D(M)\backslash e$ (for all $%
M\in \mathcal{M}$, $e\in E$).

\smallskip \ 

\noindent where the symbols ($/$) and ($\backslash $) denote the usual
operations of contraction and deletion, respectively.

\smallskip \ 

\noindent \textbf{Theorem; }Matroid duality is the unique function $D:%
\mathcal{M\longrightarrow M}$ satisfying (I)-(III).

\smallskip \ 

\noindent (see [15]-[16] for details).

That this abstract duality connects sections 2 and 5 is evident. The
difficult part is to see whether abstract duality is also related to section
3 and 4. The task seems difficult, but it may help if at least we can say
something for section 4. First, there are an indication that the Fano
matroid $F_{7}$ is connected with octonions (see Refs. [26] and [27] and
references therein). Secondly, it is known that in differential geometry the
Hodge dual is an important notion. One can show that this concept can
interpret in terms of a complete antisymmetric $\varepsilon $-symbol which
in turn one can see that this symbol admits an interpretation of a chirotope
in oriented matroid theory. One can also show that the $\varepsilon $-symbol
plays an important role in the theorem of section 4. So one wonders whether
an abstract duality definition for oriented matroid theory may be the key
notion to have a better understanding of the celebrated theorem that the
only dimensions $n$ for which we have multiplication $R^{n}\times
R^{n}\rightarrow R^{n}$ denoted with $xy=0$ implying either $x=0$ or $y=0$
are $n=1,2,4$ or $8$.

\bigskip \ 

\noindent \textbf{Acknowledgments: }This work was partially developed at the
Mathematical, Computational \& Modeling Science Center of the Arizona State
University.\ My special thanks to C. Castillo for its support.

\end{document}